\documentclass[11pt]{article}

\usepackage{doublespace}
\usepackage{epsfig}
\usepackage{amssymb}
\usepackage{verbatim}
\usepackage{delarray}
\usepackage{graphicx}
\usepackage{epic}
\usepackage{amsthm}
\usepackage{amstext}

\newcommand{\bl}{\left(}
\newcommand{\br}{\right)}
\newcommand{\re}[1]{\Re\{ {#1} \}}

\newcommand{\tr}{\mbox{Tr}}

\newcommand{\lra}{\leftrightarrow}
\newcommand{\inner}[2]{\langle{#1},{#2}\rangle}

\newcommand{\HH}{{\mathcal{H}}}
\newcommand{\R}{{\mathcal{R}}}

\newcommand{\SSS}{{\mathcal{S}}}
\newcommand{\G}{{\mathcal{G}}}

\newcommand{\F}{{\mathcal{F}}}
\newcommand{\Z}{{\mathbb{Z}}}
\newcommand{\CCC}{{\mathbb{C}}}
\newcommand{\Q}{{\mathcal{Q}}}
\newcommand{\I}{{\mathcal{I}}}

\newcommand{\ie}{{\em i.e., }}
\newcommand{\eg}{{\em e.g., }}

\newcommand{\deft}{{\stackrel{\triangle}{=}}}

\newtheorem{theorem}{Theorem}
\newtheorem{cor}{Corollary}

\newtheorem{prop}{Proposition}

\title{Geometrically Uniform Frames}
\author{Yonina C. Eldar\footnote{Research Laboratory of Electronics,
Massachusetts Institute of Technology, Room 36-615,
Cambridge, MA 02139. E-mail:
yonina@mit.edu.}\hspace*{0.03in} and 
Helmut B\"{o}lcskei\footnote{Coordinated Science Laboratory,
University of Illinois at Urbana-Champaign, Urbana, IL  61801.
Email: bolcskei@comm.csl.uiuc.edu.}}
\begin{document} 

\maketitle

%%%%%%%%%%%%%%%%%%%%%%%%%%%%%%%%%%%%%

\begin{abstract}

We introduce a new class of frames with strong symmetry properties
called {\it geometrically uniform frames (GU)}\/, that are defined over an
abelian group of unitary matrices and 
are generated by a single generating vector.
The notion of GU frames is then extended to {\it compound GU (CGU)}\/ frames 
which are generated 
by an  abelian group of unitary matrices using multiple
generating vectors.

The dual frame vectors and canonical tight
frame vectors associated with GU frames are shown to be GU
and  therefore generated by a single generating vector, which 
can be computed very efficiently using a Fourier
transform defined over the generating group of the frame. Similarly,
the dual frame vectors and canonical tight
frame vectors associated with CGU frames are shown to be CGU.

The impact of removing single or multiple elements from a GU frame is
considered. 
A  systematic method for constructing optimal GU frames
from a given set of frame vectors that are not GU is also developed.
Finally, the Euclidean distance properties of GU frames
are discussed and conditions are derived on the abelian group of
unitary matrices to yield 
GU frames with strictly positive distance spectrum irrespective of
the generating vector. 
\end{abstract}

%%%%%%%%%%%%%%%%%%%%%%%%%%%%%%%%%%%%

%%%%%%%%%%%%%%%%%%%%%%%%%%%%%%%
\section{Introduction}
%%%%%%%%%%%%%%%%%%%%%%%%%%%%%%%%

Frames are generalizations of bases which lead to redundant signal
expansions \cite{DS52,Y80}.  A (finite) frame for a Hilbert space
$\HH$ is a set of not necessarily linearly independent vectors that spans 
$\HH$.
Since the frame vectors can be linearly
dependent, the conditions on frame vectors are usually not as stringent
as the conditions on bases, allowing for increased flexibility in
their design \cite{HW89,A95}.

Frames were first introduced by Duffin and Schaeffer \cite{DS52} in
the context of nonharmonic Fourier series, and play an important role
in the theory of nonuniform sampling \cite{DS52,Y80,B92} and wavelet
theory \cite{HW89,D90}. Recently,  frames
have been used to analyze and design oversampled filter banks
\cite{B97,BHF98,Cve98} and error correction codes \cite{F99}. Frames
have also been applied to the development of modern uniform and
nonuniform sampling techniques \cite{AG01}, to various detection
problems \cite{EO01,EC01}, and to multiple description source coding 
\cite{GKV98}.

Two important classes of highly structured frames are Gabor (Weyl-Heisenberg 
(WH))
frames \cite{G46,fei_book} and 
wavelet frames \cite{HW89,D90,D92}. Both classes of frames are generated
by a single generating function. WH frames are obtained by translations
and modulations of the generating function (referred to as the window
function), and wavelet frames are obtained 
by shifts and dilations of the generating function (referred to as the
mother wavelet).
In Section~\ref{sec:guframes} of this paper, we introduce a new class of frames
which we refer to 
as geometrically uniform (GU) frames, that like WH and wavelet frames are
generated from a single generating vector. These frames are defined by
an abelian group $\Q$ of unitary matrices, referred to as the generating
group of the frame. We note that WH frames and wavelet frames are in general
not GU since the underlying group of matrices is in general not abelian.
GU frames are based on the notion of GU vector sets first introduced by Slepian \cite{S68}
and later extended by Forney \cite{F91}, which
are known to  have strong
symmetry properties that may be desirable in various
applications such as channel coding \cite{F91,U82,R98}.

The notion of GU frames is then extended to frames that are generated
by an abelian group $\Q$ of unitary matrices using {\it multiple}\/
generating vectors.  Such frames are not necessarily GU, but consist of
subsets of GU vector sets that are each generated by
$\Q$. We refer to this class of frames as compound GU (CGU) frames, 
and develop their properties in Section~\ref{sec:cguf}. CGU frames are a generalization
of filter bank frames introduced in \cite{B97,BHF98,Cve98}.
An interesting class of
frames results when the set of generating vectors is
itself GU, generated by an abelian group $\G$. 
(Note that this class of frames will in general not be GU).
As we show, these frames
are a generalization of WH frames in which $\Q$ is the group of
translations and $\G$ is the group of 
modulations.

Given a frame for $\HH$, any signal in $\HH$ can be represented as a
linear combination of the frame vectors. However, if the frame vectors
are linearly dependent, 
then the coefficients in this expansion are not unique. A popular
choice of coefficients are the inner products of the signal with a set of
analysis frame vectors  called the dual frame vectors \cite{D92}. This
choice of coefficients has
the property that among all possible coefficients it has the minimal
$l_{2}$-norm \cite{D92,K94}. 

In Section~\ref{sec:dual}, we show that the dual frame vectors
associated with a GU frame are also GU, 
and  therefore generated by a single generating vector. Furthermore, we 
demonstrate that
the generating vector can be computed very efficiently using a Fourier
transform defined over the generating group $\Q$ of the frame.
Similarly, in Section~\ref{sec:cguf} we show that the dual frame
vectors associated with a CGU frame  are also CGU. 
When the generating vectors of the CGU frame are GU and generated by a
group $\G$ that commutes up to a phase factor with the group $\Q$, 
the dual frame is generated by a
{\em single} generating vector, a result well known for WH frames.

An important topic in frame theory is the behavior of a frame when
elements of the frame are removed. In
Section~\ref{sec:pruning}, we show that the
frame bounds of the frame resulting from removing a single vector
of a GU frame are the same regardless of the  particular vector
removed. In this sense GU frames exhibit an interesting robustness
property which is of particular importance in applications such as
multiple description source coding \cite{GKV98}. We also
consider the behavior of a GU frame when groups of frame
elements are removed.

In the special case of a 
tight frame the dual frame vectors are proportional to the original
frame vectors so that the reconstruction 
formula is particularly simple. In many applications it is 
therefore  desirable to construct a tight
frame from an arbitrary set of frame vectors.
A popular tight frame construction 
is the so-called canonical tight frame \cite{D92,B97,BJ99,JB00,EF01,JS00},
first proposed in the context of wavelets in \cite{M86}.  The
canonical tight frame is relatively simple to construct, it is optimal in a
least-squares sense \cite{EF01f,JS00,S00b}, it can be
determined directly from the given vectors, and plays an important
role in wavelet theory \cite{D88,M89,UA93}.  
Like the  dual frame vectors, we show that the canonical tight frame
vectors associated with a GU frame are GU, and the 
canonical tight frame vectors associated with a CGU frame are 
CGU. When the generating vectors of the CGU frame are GU and generated by a
group $\G$ that commutes up to a phase factor with  $\Q$, 
the canonical tight frame vectors can be obtained by a single generating
vector, generalizing a result well-known in WH frame theory.

Since GU frames have nice symmetry properties, it may be desirable to
construct such a frame from a given set of frame vectors.
The problem of frame design has received relatively little attention
in the literature. Systematic methods for constructing
optimal tight frames have been considered
\cite{EF01f,JS00,S00b}. Methods for generating frames starting from a
given 
frame are described in \cite{A95}. 
In Section~\ref{sec:optimal}, we systematically construct optimal
GU frames from a given set of vectors, that are closest in a
least-squares sense to the original frame vectors.
We consider three different constraints on the
GU frame vectors. First, we treat the case in which the inner
products of the frame vectors are known.
The optimizing frame is referred to as the
scaled-constrained  least-squares GU frame (SC-LSGUF). Next, we
consider the case where the inner products are known up to a scale
factor. The  optimizing frame in this case is referred to as the
constrained  least-squares GU frame (C-LSGUF). Finally, we consider the
case in which both the inner products and the scaling are chosen to
minimize the least-squares error between the original frame and the resulting 
tight frame.
The optimizing frame is the least-squares GU
frame (LSGUF). 

In Section~\ref{sec:distance}, we consider distance properties of 
GU frames, which may be of interest
when using GU frames for code design (group codes) \cite{S68,F91}. In
particular, we introduce a class of GU 
frames with strictly positive distance spectra
for all choices of generating vectors. Such GU frames are shown to be 
generated
by fixed-point-free groups \cite{S00}.

Before proceeding to the detailed development, in
Section~\ref{sec:frames} we provide a brief  introduction to frame
expansions.

%%%%%%%%%%%%%%%%%%%%%%%%%%%%%%%%%%%%%%%%%%%%%%%%%
\section{Frames}
\label{sec:frames}
%%%%%%%%%%%%%%%%%%%%%%%%%%%%%%%%%%%%%%%%%%%%%%%%

Frames, which are generalizations of bases, were introduced in the
context of nonharmonic Fourier series by Duffin and Schaeffer
\cite{DS52} (see also \cite{Y80}). Recently, the theory of frames has
been expanded \cite{HW89,D90,D92,A95}, in part due to the utility of
frames in analyzing wavelet decompositions.

Let $\{\phi_i,1 \leq i \leq n \}$ denote a set of $n$ complex vectors in an
$m$-dimensional Hilbert space $\HH$.  The vectors
$\phi_i$ form a {\em frame} for $\HH$ if there exist constants $A > 0$
and $B < \infty$ such that
\begin{equation}
\label{eq:frame}
A||x||^2 \leq \sum_{i=1}^n |\inner{x}{\phi_i}|^2 \leq B ||x||^2,
\end{equation}
for all $x \in \HH$ \cite{D92}. 
In this paper, we restrict our attention to the case where $m$ and $n$ are 
finite.
The lower bound in (\ref{eq:frame}) ensures that the vectors $\phi_i$
span $\HH$; thus we must have $n \geq m$. Since $n<\infty$ the right
hand inequality of (\ref{eq:frame}) is always satisfied with
$B=\sum_{i=1}^n ||\phi_{i}||^2$, so that any finite set of
vectors that spans $\HH$ is a frame for $\HH$. In particular, any basis
for $\HH$ is a frame for $\HH$. However, in contrast to basis vectors
which are linearly independent, frame vectors with $n>m$ are linearly
dependent.  If the bounds $A=B$ in (\ref{eq:frame}), then the frame is
called a {\it tight frame}\/. If in addition $A=B=1$, then  the frame
is called  a
{\em normalized tight frame}.  The redundancy of the frame is defined as
$r=n/m$, \ie $n$ vectors in an $m$-dimensional space.

The {\it frame operator} corresponding to the frame vectors
$\{\phi_i,1 \leq i \leq n\}$ is defined
as \cite{D92}
\begin{equation}
\label{eq:S}
S=\sum_{i=1}^n \phi_i\phi_i^*=\Phi\Phi^*,
\end{equation}
where $\Phi$ is the matrix of columns $\phi_i$, and $(\cdot)^*$
denotes the Hermitian transpose. Using the frame operator, (\ref{eq:frame})
can be rewritten as
\begin{equation}
\label{eq:frames}
A||x||^2 \leq \inner{Sx}{x} \leq B ||x||^2.
\end{equation}

From (\ref{eq:frames}) it follows that the  tightest possible frame
bounds $A$ and $B$ are given by 
$A=\min_i \lambda_i(S)$ and $B=\max_i
\lambda_i(S)$, where $\{\lambda_i(S),1 \leq i \leq m\}$ are the
eigenvalues of the frame operator $S$. Throughout the paper,
when referring to ``the frame bounds'' we implicitly assume
the tightest possible frame bounds unless otherwise stated.

If the vectors $\{\phi_i,1 \leq i \leq n\}$ form a frame for $\HH$, then
any $x \in \HH$ can be expressed as a linear combination of these
vectors: $x=\sum_{i=1}^n a_i \phi_i$. If $n>m$ the coefficients in this
expansion are not unique. A possible choice is\footnote{We use the
notation $\inner{x}{y}=x^*y$.} 
$a_i=\inner{\bar{\phi}_i}{x}$  where $\bar{\phi}_i$ are the {\em dual
frame vectors} \cite{D92} of the frame vectors $\phi_i$, and are given
by 
\begin{equation}
\label{eq:dual}
\bar{\phi}_i=S^{-1}\phi_i.
\end{equation}
We can express $S^{-1}$ directly
in terms of $S$ as \cite{D90} 
\begin{equation}
\label{eq:sinv}
S^{-1}=\frac{2}{A+B}\left[I-\left(I-\frac{2}{A+B}S
\right)\right]^{-1}=\frac{2}{A+B}\sum_{l=0}^\infty \bl
I-\frac{2}{A+B}S \br^l, 
\end{equation}
where $A$ and $B$ are the frame bounds of 
$\{\phi_i\}$.
The choice of coefficients $a_i=\inner{\bar{\phi}_i}{x}$ has
the property that among all possible coefficients it has the minimal 
$l_{2}$-norm
\cite{D92,K94}. 

There are other choices of dual frame vectors $y_i$ such that for any $x
\in \HH$, $x=\sum_{i=1}^n \inner{y_i}{x} \phi_i$. Specifically, with
$Y$ denoting 
the matrix of columns $y_i$, any other choice corresponds to $Y$ of
the form \cite{BG80} 
\begin{equation}
Y=S^{-1}\Phi+T \bl I-\Phi^* S^{-1}\Phi \br,
\end{equation}
where $T$ is an arbitrary matrix with bounded elements.
However, the particular choice $y_i=\bar{\phi}_i$ has some desirable
properties. 
Besides resulting in the minimal $l_{2}$-norm coefficients, in many cases
the choice $y_i=\bar{\phi}_i$ yields  
frame vectors that share the same symmetries as the original frame
vectors. Specifically, in Section~\ref{sec:dual} we show that the dual
frame vectors associated with a geometrically uniform (GU) frame  are
also GU, and in Section~\ref{sec:cguf} we show that 
the dual
frame vectors associated with a compound GU (CGU) frame are also
CGU. Finally,
in the case of a tight frame the dual frame vectors lead to a
particularly simple expansion. Specifically, in this case
$S=AI_m$ so that
$S^{-1}=(1/A)I_m$, and the 
dual frame vectors are  $\{\bar{\phi}_i=(1/A)\phi_i,1 \leq i \leq
n\}$. Since a tight frame expansion of a signal is very simple,
it is popular in many applications \cite{D92}.

Suppose we are given a set of vectors $\{\phi_i,1 \leq i \leq n\}$
that form a frame 
for $\HH$, with frame bounds $A \neq B$. It may then be desirable to
construct a tight frame from these vectors.
A popular tight frame construction 
is the canonical tight frame \cite{D92,B97,BJ99,JB00,JS00,EF01f},
first proposed in the context of wavelets in \cite{M86}.  
The {\em canonical tight frame vectors} $\{\mu_i,1 \leq i \leq n\}$
associated with the 
vectors $\{\phi_i,1 \leq i \leq n\}$ are given by
\begin{equation}
\label{eq:canonical}
\mu_i=S^{-1/2}\phi_i,
\end{equation}
where $S^{-1/2}$ is the positive-definite square root of $S^{-1}$. 
We can express $S^{-1/2}$ directly
in terms of $S$ as \cite{daubjj91}
\begin{equation}
\label{eq:ssqrt}
S^{-1/2}=\sqrt{\frac{2}{A+B}}\sum_{l=0}^\infty
\frac{(2l)!}{2^{2l}(l!)^2}\bl I-\frac{2}{A+B}S \br^l, 
\end{equation}
where $A$ and $B$ are the frame bounds associated with the frame
$\{\phi_i\}$. Note,
that $\{US^{-1/2}\phi_i\}$ with $U$ an arbitrary unitary matrix yields a tight
frame as well.
The canonical tight frame, however, has the property that it is the 
closest normalized tight frame to the vectors $\{\phi_i\}$ in a
least-squares sense \cite{EF01,JS00,S00b}. 

From (\ref{eq:dual})  and (\ref{eq:canonical}) we see that in
order to compute the dual frame vectors and the canonical tight frame vectors
associated with a  frame $\{\phi_i\}$, we need to compute the
matrices $S^{-1}$ and $S^{-1/2}$ and then apply them to each of the frame
vectors $\phi_i$. In the next section, we introduce a class
of frames  
that have strong symmetry properties called
{\em geometrically uniform (GU) frames}.
As we show in Section~\ref{sec:dual}, the dual frame vectors and the
canonical tight frame vectors 
associated with a GU frame are generated by a single generating function,
and can therefore be computed very
efficiently.

%%%%%%%%%%%%%%%%%%%%%%%%%%%%%%%%%%%%%%%%%%%%%%%%%
\section{Geometrically Uniform Frames}
\label{sec:guframes}
%%%%%%%%%%%%%%%%%%%%%%%%%%%%%%%%%%%%%%%%%%%%%%%%

A set of vectors $\SSS = \{\phi_i, 1 \leq i \leq n\}$ is geometrically
uniform (GU) \cite{F91,S68,EF01} if every vector in the set has the
form $\phi_i=U_i\phi$, 
where $\phi$ is an arbitrary {\em generating vector} and the matrices
$\{U_i,1 \leq i 
\leq n\}$ are unitary and form an abelian group\footnote{That is, $\Q$
contains the identity matrix $I$; if $\Q$ contains $U_i$, then it also
contains its inverse $U_i^{-1}$; the product $U_i U_j$ of any
two elements of $\Q$ is again in $\Q$; and $U_i U_j=U_j U_i$ for any two
elements in $\Q$ \cite{A88}.}  $\Q$. For concreteness we assume that
$U_1=I$ so that $\phi_1=\phi$.
The group $\Q$ will be called the \emph{generating group} of
$\SSS$.

Alternatively, a vector set is GU if given any two vectors $\phi_i$
and $\phi_j$  
in the set, there is an
isometry (a norm-preserving linear transformation) $Z_{ij}$
that transforms $\phi_i$ into $\phi_j$ while leaving the set
invariant \cite{F91}. Thus for every $i$, $Z_{ij}\phi_i=\phi_j$.
Intuitively, a vector set is GU if it
``looks the same'' geometrically from any of the points in the set. 
Some examples of GU vector sets are considered in \cite{F91}.

A set of vectors $\{\phi_i \in \HH,1 \leq i \leq n\}$ forms a
{\em geometrically uniform frame} 
for $\HH$, if the vectors $\{\phi_i\}$ are GU and span $\HH$.

As we show in the following proposition, the frame bounds of a GU
frame can be bounded by
the norm of the generating vector. 
\begin{prop}
\label{prop:bounds}
Let $\SSS = \{\phi_i = U_i\phi, U_i \in \Q\}$ be a geometrically
uniform frame with frame bounds $A$ and $B$, where $\phi$ is an
arbitrary generating vector. Then
$A\leq \frac{n}{m} ||\phi||^2 \leq B$. 
If in addition the frame is tight, then
$A=B=\frac{n}{m} ||\phi||^2 $. 
\end{prop}
\begin{proof}
We can express the frame operator corresponding to the frame vectors
$\phi_i$ as
\begin{equation}
S=\sum_{i=1}^n U_i\phi\phi^*U_i^*.
\end{equation}
Then,
\begin{equation}
\tr(S)=\sum_{i=1}^n \tr(U_i\phi\phi^*U_i^*)=n  ||\phi||^2,
\end{equation}
so that
\begin{equation}
\sum_{i=1}^m \lambda_i(S)=\tr(S)=n ||\phi||^2.
\end{equation}
Therefore,
\begin{equation}
A=\min_i \lambda_i(S) \leq \frac{1}{m}\sum_{i=1}^m
\lambda_i(S)=\frac{n}{m}||\phi||^2,
\end{equation}
and
\begin{equation}
B=\max_i \lambda_i(S) \geq \frac{1}{m}\sum_{i=1}^m
\lambda_i(S)=\frac{n}{m}||\phi||^2. 
\end{equation}
\end{proof}

Since $U_i^* = U_i^{-1}$, the inner product of two  vectors in $\SSS$ is
\begin{equation}
\inner{\phi_i}{\phi_j} = \phi^* U_i^{-1} U_j\phi=
s(U_i^{-1}U_j),
\end{equation}
where $s$ is the function on $\Q$ defined by
\begin{equation}
s(U_i) = \phi^*U_i \phi.
\end{equation}
For fixed $i$, the set $U_i^{-1}\Q = \{U_i^{-1}U_j, U_j \in \Q\}$ 
is just a permutation of $\Q$ since $U_i^{-1}U_j \in \Q$
for all $i,j$ \cite{A88}.  Therefore,
the $n$ numbers $\{s(U_i^{-1}U_j), 1 \leq j \leq n\}$ are a permutation
of the  numbers $\{s(U_i), 1 \leq i \leq n\}$. The same is true for fixed
$j$. Consequently, every row and column of 
the $n \times n$ Gram matrix $G =
\{\inner{\phi_i}{\phi_j}\}$ 
is a permutation of the numbers $\{a_i=s(U_i), 1 \le i \le n\}$. 

A matrix $G$ whose rows (columns) are a permutation of the first row
(column) will be called a
{\em permuted matrix}\footnote{An example of a permuted
matrix is
\begin{equation}
\label{eq:guex}
\left[
\begin{array}{rrrr}
a_1 & a_2& a_3 &a_4 \\
a_2& a_1& a_4& a_3\\
a_3& a_4& a_1& a_2\\
a_4& a_3& a_2& a_1
\end{array}
\right].
\end{equation}}.
Thus, we have shown that the Gram matrix of a GU
vector set is a permuted matrix. 
Furthermore, if the Gram matrix $G=\{\inner{\phi_i}{\phi_j}\}$ is a permuted 
matrix
and in addition $G=G^{T}$,  then the vectors $\{\phi_i\}$ are GU
\cite{E01}. We therefore have the following proposition.
\begin{prop}
\label{prop:permuted}
The Gram matrix $G=\{\inner{\phi_i}{\phi_j}\}$ corresponding to a
geometrically uniform
vector set  
$\SSS = \{\phi_i \in \HH, 1 \leq i \leq n\}$ is a permuted matrix. Conversely,
if the  Gram matrix  $G=\{\inner{\phi_i}{\phi_j}\}$ is
a permuted matrix, and $\inner{\phi_i}{\phi_j}=\inner{\phi_j}{\phi_i}$
for all $i,j$, then the vectors $\{\phi_i\}$ are geometrically
uniform. If in addition  
the vectors $\{\phi_i\}$ span $\HH$, then they form a geometrically
uniform frame for $\HH$.
\end{prop}

It will be convenient to replace the multiplicative group $\Q$ by an
additive group $Q$ to which $\Q$ is isomorphic\footnote{ Two groups
$\Q$ and $\Q'$ are 
{\em isomorphic}, denoted by $\Q \cong \Q'$, if there is a bijection
(one-to-one and onto map) $\varphi: \Q \to \Q'$ which satisfies
$\varphi(xy)=\varphi(x)\varphi(y)$ for all $x,y \in \Q$ \cite{A88}.}.
Specifically, it is well known (see \eg \cite{A88}) that
every finite abelian group $\Q$ is isomorphic to a direct product $Q$
of a finite number of cyclic groups: $\Q \cong Q = \Z_{n_1} \times
\cdots \times \Z_{n_p}$, where $\Z_{n_t}$ is the cyclic additive group
of integers modulo $n_t$, and $n = \prod_t n_t$.  Thus every element
$U_i \in \Q$ can be associated with an element $q \in Q$ of the form
$q = (q_1,q_2,\ldots ,q_p)$, where $q_t \in \Z_{n_t}$; this
correspondence is denoted by $U_i \lra q$.

Each vector $\phi_i=U_i\phi$ is then denoted as $\phi(q)$, where $U_i \lra
q$.
The zero element $0=(0,0,\ldots,0) \in Q$ corresponds to the identity
matrix $I \in \Q$, and an additive inverse $-q \in Q$ corresponds to a
multiplicative inverse $U_i^{-1} = U_i^* \in \Q$.
The Gram matrix is then the $n \times n$ matrix 
\begin{equation}
\label{eq:sij}
G = \{\inner{\phi(q')}{\phi(q)}, q', q \in Q\} = \{s(q - q'),
q', q \in Q\},
\end{equation}
with row and column indices $q', q \in Q$, where $s$ is now the function
on $Q$ defined by
\begin{equation}
s(q) = \inner{\phi(0)}{\phi(q)}.
\end{equation}

The Fourier transform (FT) of a complex-valued function $\varphi: Q
\to \CCC$ defined on $Q = \Z_{n_1} \times \cdots \times \Z_{n_p}$ is the
complex-valued function $\hat{\varphi}: Q \to \CCC$ defined by
\begin{equation}
\label{eq:fh}
\hat{\varphi}(h) = \frac{1}{\sqrt{n}}\sum_{q \in Q} \inner{h}{q}
\varphi(q),
\end{equation}
where the Fourier kernel $\inner{h}{q}$ is
\begin{equation}
\label{eq:hq}
\inner{h}{q} = \prod_{t=1}^p e^{-2 \pi i h_tq_t/n_t}.
\end{equation}
Here $h_t$ and $q_t$ are the $k$th components of $h$ and $q$
respectively, and the product $h_tq_t$ is taken as an ordinary integer
modulo $n_t$.

The FT matrix over $Q$ is defined as the $n \times n$ matrix $\F =
\{\frac{1}{\sqrt{n}}\inner{h}{q}, h,q \in Q\}$.  The FT of a column
vector $\varphi = \{\varphi(q), q \in Q\}$ is then the column vector
$\hat{\varphi} = \{\hat{\varphi}(h), h \in Q\}$ given by
$\hat{\varphi} = \F\varphi$.  Since $\F$ is unitary,
we obtain the inverse FT formula
\begin{equation}
\varphi = \F^*\hat{\varphi} =
\left\{\frac{1}{\sqrt{n}}\sum_{h \in Q}
\inner{h}{q}^* \hat{\varphi}(h),q \in Q\right\}.
\end{equation}

\newpage 
As we show in the following theorem, the FT matrix plays an important
role in defining GU frames.
\begin{theorem}
\label{thm:guft}
A set of vectors $\{\phi_i,1 \leq i \leq n\}$ in an $m$-dimensional
Hilbert space $\HH$ is geometrically
uniform if and only if the Gram matrix $G=\{\inner{\phi_i}{\phi_j}\}$
is diagonalized by a Fourier transform matrix $\F$ over a finite
product of cyclic groups $Q$. The
vectors $\{\phi_i\}$ form a geometrically uniform frame for $\HH$ if
 in addition $G$ has rank $m$.
\end{theorem}
\begin{proof}
The vectors $\{\phi_i\}$ form a frame for $\HH$ if and only if they
span $\HH$, which implies that the rank of $G$ must be equal to
$m$. 

For a  GU vector set with generating group $\Q \cong
Q$, the FT over $Q$ 
diagonalizes the Gram matrix $G$ \cite{EF01}. Thus, to complete the proof of
the theorem 
we  need to prove that if $G$ is diagonalized by a FT matrix
$\F$ over the group $Q$, then the vector set $\{\phi(q),q \in Q\}$ is GU.

Let $\Phi$ be the matrix of columns $\phi(q)$, so that $G=\Phi^*\Phi$. 
Since $\F$ diagonalizes $G$, $G$ has an eigendecomposition of the form
$G=\F D\F^*$ for a diagonal matrix $D$ with  diagonal
elements $d_i$, where the first $m$ diagonal elements may be non-zero and the
remaining diagonal elements are all zero.
Then $\Phi$ has an SVD of
the form $\Phi=U\Sigma \F^*$ where $U$ is an arbitrary unitary matrix and
$\Sigma$ is an $m \times n$ diagonal matrix with diagonal elements
$\{\sigma_i=\sqrt{d_i},1 \leq i \leq m\}$.

Let $\{f(q),q \in Q\}$ denote the columns of $\F^*$. From the
definition of $\F$, the components of $f(0)$ are all equal
$1/\sqrt{n}$, and
$f(q)=B(q) f(0)$ where $B(q)$ is a diagonal unitary matrix with
diagonal elements $\{\inner{h}{q},h \in Q\}$, where $\inner{h}{q}$
is given by (\ref{eq:hq}). Then,
\begin{equation}
\phi(q)=U\Sigma f(q)=U\Sigma B(q)f(0)=UB(q) \Sigma f(0)=UB(q)U^* \phi,
\end{equation}
where $\phi=U\Sigma f(0)$, and where we used the fact that diagonal matrices
commute. If we now define $U(q)=UB(q)U^*$, then we have that
$\phi(q)=U(q)\phi$ where the matrices $\{U(q),q \in Q\}$ are
unitary. 

We now show that the group $\Q=\{U(q),q \in Q\}$ is an
abelian group.
First, we have that $U(0)=UB(0)U^*=UU^*=I$ so that $I \in \Q$. Next,
$U^{-1}(q)=UB^{-1}(q)U^*=UB(-q)U^*$ so that $U^{-1}(q) \in \Q$ since 
$-q \in Q$. Finally, $U(q)U(h)=UB(q)B(h)U^*=UB(h)B(q)U^*=U(h)U(q)$
since diagonal matrices commute, and $U(q)U(h) \in \Q$ since
$B(h)B(q)=B(h+q)=B(q')$ for some $q' \in Q$.

We therefore conclude that $\phi(q)=U(q)\phi$ where the matrices $U(q)$
are unitary and 
form an abelian group, so that the vectors $\phi(q)$ are geometrically
uniform.
\end{proof}

As a consequence of Theorem~\ref{thm:guft} we have the following
corollary.
\begin{cor}
\label{cor:gusvd}
A set of vectors $\{\phi_i \in \HH, 1 \leq i \leq n\}$ is geometrically
uniform if and only if  
the matrix $\Phi$ of columns $\phi_i$ has an SVD of the form
$\Phi=U\Sigma \F^*$, where $U$ is an arbitrary unitary matrix, $\Sigma$
is an arbitrary diagonal matrix with diagonal elements $\{\sigma_i,1 
\leq i \leq m\}$, and $\F$ is a FT matrix over a
direct product of cyclic groups. In addition, the vectors $\{\phi_i\}$
form a geometrically uniform frame for $\HH$ if they span
$\HH$ or, equivalently, if $\sigma_i \neq 0$ for $1 \leq i \leq m$.
\end{cor}

%%%%%%%%%%%%%%%%%%%%%%%%%%%%%%%%%%%%%%%%%%%%%%%%%
\section{Dual and Canonical Tight Frames Associated With GU Frames}
\label{sec:dual}
%%%%%%%%%%%%%%%%%%%%%%%%%%%%%%%%%%%%%%%%%%%%%%%%

In  Section~\ref{sec:construction}, we show that  the
dual frame vectors and the canonical tight frame vectors associated with a
GU frame are also GU. This property can then be used to compute
the dual and canonical tight frames very efficiently. Further properties of
the canonical tight frame vectors are discussed in
Section~\ref{sec:properties}.

%%%%%%%%%%%%%%%%%%%%%%%%%%%%%%%%%%%%%%%%%%%%%%%%%
\subsection{Constructing the dual and canonical tight frames}
\label{sec:construction}
%%%%%%%%%%%%%%%%%%%%%%%%%%%%%%%%%%%%%%%%%%%%%%%%

Let $\SSS = \{\phi_i = U_i\phi, U_i \in \Q\}$ be a GU 
frame generated by a finite abelian group $\Q$ of unitary
matrices, where $\phi$ is an arbitrary generating vector.
Then the
frame operator $S$ defined by (\ref{eq:S}) commutes with each of the
unitary matrices $U_i$ in 
the generating group $\Q$. Indeed,
expressing the frame operator as
\begin{equation}
S=\sum_{i=1}^n U_i\phi\phi^*U_i^*,
\end{equation}
we have that for all $j$,
\begin{eqnarray}
SU_j & = & \sum_{i=1}^n U_i\phi_i\phi_i^*U_i^*U_j \nonumber \\
& = & U_j \sum_{i=1}^n U_j ^*U_i\phi\phi^*U_i^*U_j \nonumber \\
& = & U_j \sum_{i=1}^n U_i\phi\phi^*U_i \nonumber \\
& = & U_jS,
\end{eqnarray}
since $\{U_j ^*U_i,1 \leq i \leq n\}$ is just a permutation of $\Q$.

If $S$ commutes with $U_j$, then from (\ref{eq:sinv}) and
(\ref{eq:ssqrt}) we have that  $S^{-1}$ and $S^{-1/2}$ also commute
with $U_j$ for all $j$. Thus,
\begin{equation}
\bar{\phi}_i=S^{-1}\phi_i=S^{-1}U_i\phi=U_iS^{-1}\phi=U_i \bar{\phi},
\end{equation}
where $\bar{\phi}=S^{-1} \phi$, which shows that the dual frame vectors
$\{\bar{\phi}_i=S^{-1}\phi_i\}$ are GU with generating group equal
to $\Q$. 

Similarly,
\begin{equation}
\mu_i=S^{-1/2}\phi_i=S^{-1/2}U_i\phi=U_iS^{-1/2}\phi=U_i \mu,
\end{equation}
where $\mu=S^{-1/2} \phi$, which shows that the canonical tight frame vectors 
$\{\mu_i=S^{-1/2}\phi_i\}$ are also GU with generating group $\Q$.

Therefore, to compute the dual frame vectors or the canonical tight frame
vectors all we
need is to compute  the generating vectors $\bar{\phi}$ and $\mu$,
respectively. 
The remaining frame vectors are then obtained by applying the group $\Q$ to the
corresponding generating vectors. 

We now show that the generating vectors can be computed very
efficiently using the FT. From Corollary~\ref{cor:gusvd} we
have that $\Phi$
has an SVD of the form
\begin{equation}
\label{eq:svdphi}
\Phi = U \Sigma \F^* = \sum_{h \in Q}
\sigma(h)u(h)\F^*(h).
\end{equation}
Here $\Sigma$ is a diagonal matrix with
diagonal elements $\{\sigma(h) = n^{1/4}\sqrt{\hat{s}(h)}, h \in Q\}$
where $\{\hat{s}(h), h \in Q\}$ is the FT of $\{s(q), q \in Q\}$,
$U$ is the 
matrix of columns $u(h)$, where
\begin{eqnarray}
\label{eq:uh}
u(h) & = & \left\{\begin{array}{ll}
\Phi \F(h)/\sigma(h) = \hat{\phi}(h)/\sigma(h),  &
 \mbox{if}\,\, \sigma(h) \neq 0; \\
0, & \mathrm{otherwise}
\end{array}\right.
\end{eqnarray}
with
\begin{equation}
\label{eq:ftphi}
\hat{\phi}(h) = \frac{1}{\sqrt{n}}\sum_{q \in Q}
\inner{h}{q}\phi(q)
\end{equation}
denoting the $h$th element of the FT of $\Phi$ regarded as a
row vector of column vectors, \mbox{$\Phi = \{\phi(q), q \in G\}$},
and $\F^* =
\{\frac{1}{\sqrt{n}}\inner{h}{q}^*, h, q \in Q\}$ has rows $\F^*(h) =
\{\frac{1}{\sqrt{n}}\inner{h}{q}^*, q \in Q\}$.  

It then follows that
\begin{equation}
\label{eq:barphi}
\bar{\phi}=S^{-1}\phi=\frac{1}{\sqrt{n}}\sum_{h \in \I}
\frac{1}{\sigma(h)}u(h), 
\end{equation}
where $h \in \I$ if $\sigma(h) \neq 0$. Similarly,
\begin{equation}
\label{eq:mu}
\mu=S^{-1/2}\phi=\frac{1}{\sqrt{n}}\sum_{h \in \I} u(h).
\end{equation}

\newpage
We summarize our results in the following theorem:
\begin{theorem}[GU frames]
\label{thm:gu}
Let $\SSS = \{\phi_i = U_i\phi, U_i \in \Q\}$ be a GU 
frame generated by a finite abelian group $\Q$ of unitary
matrices, where $\phi$ is an arbitrary generating vector, and let
$\Phi$ be the matrix 
of
columns $\phi_i$. Let $Q$ be an additive abelian group isomorphic to $\Q$,
let
$\{\phi(q), q \in Q\}$ be the elements of $\SSS$ under this isomorphism, and
let 
$\F$ be the Fourier transform matrix over $G$.
Then 
\begin{enumerate}
\item the dual frame vectors $\{\bar{\phi}_i,1 \leq i \leq n\}$ are
geometrically uniform with generating group $\Q$ and generating vector
$\bar{\phi}=(1/\sqrt{n})\sum_{h \in \I} (1/\sigma(h))u(h)$, where
\begin{enumerate}
\item $\{\sigma(h) = n^{1/4}\sqrt{\hat{s}(h)}, h \in Q\}$ are
the singular values of $\Phi$,
\item $\{\hat{s}(h), h \in Q\}$ is the Fourier transform of the
inner-product
sequence \mbox{$\{\inner{\phi(0)}{\phi(q)}, q \in Q\}$},
\item $\I$ is the set of indices $h \in Q$ for which $\sigma(h) \neq 0$,
\item $u(h) = \hat{\phi}(h)/\sigma(h)$ for $h \in \I$,
\item $\{\hat{\phi}(h), h \in
Q\}$ is the Fourier transform of $\{\phi(q), q \in Q\}$,
\end{enumerate}
\item the canonical tight frame vectors $\{\mu_i,1 \leq i \leq n\}$ are
GU with generating group $\Q$ and generating vector
$\mu=(1/\sqrt{n})\sum_{h \in \I} u(h)$;
\item the frame bounds of the frame $\{\phi_i,1 \leq i \leq n\}$ are given by
$A=\sqrt{n}\min_{h \in \I}\hat{s}(h)$ and $B=\sqrt{n}\max_{h 
\in \I}\hat{s}(h)$. 
\end{enumerate}
\end{theorem}

An important special case of Theorem~\ref{thm:gu} 
is the case in which
the  generating group $\Q$ is {\em cyclic} so that $U_i = V^{i-1},
1 \le i \le n$, where $V$ is a unitary matrix with $V^n = I$.  A cyclic
group generates a cyclic vector set $\SSS =
\{\phi_i=V^{i-1}\phi,\,\,1
\leq i \leq n\}$, where $\phi$ is arbitrary.
If $\Q$ is cyclic, then $G$ is a circulant
matrix\footnote{A circulant matrix is a  matrix where every row (or
column) is obtained by a right circular shift (by one position) of the
previous row (or column). An example is:
$\left[ \begin{array}{ccc}
a_0 & a_2 & a_1 \\
a_1 & a_0 & a_2 \\
a_2 & a_1 & a_0
\end{array} \right].$}, and
$Q$ is the cyclic group $\Z_n$. The FT kernel is then $\inner{h}{g} =
e^{-2 \pi ihg/n}$ for $h,g \in \Z_n$, and the FT matrix $\F$ reduces
to the $n \times n$ DFT matrix.  The singular values of $\Phi$ are
then $n^{1/4}$ times the square roots of the DFT values of the inner
products $\{\inner{\phi_1}{\phi_j},1 \leq j \leq n\}$.

%%%%%%%%%%%%%%%%%%%%%%%%%%%%%%%%%%%%%%%%
\subsection{Properties of the canonical tight frame}
\label{sec:properties}
%%%%%%%%%%%%%%%%%%%%%%%%%%%%%%%%%%%%%%

The canonical tight frame vectors $\mu_i$ corresponding to the frame vectors
$\phi_i$ have the property that they are the closest normalized tight
frame vectors 
to the vectors $\phi_i$, in a least-squares sense \cite{EF01f,JS00,S00b}. Thus,
the vectors $\mu_i$ 
are the normalized tight frame vectors that minimize the least-squares error 
\begin{equation}
\sum_{i=1}^n \inner{\phi_i-\mu_i}{\phi_i-\mu_i}.
\end{equation}

We now show that when the original frame vectors $\phi_i$ are
GU with generating group $\Q$, the canonical tight frame
vectors have the additional 
property that from all normalized tight frame vectors they maximize 
\begin{equation}
\label{eq:rpmu}
R_{\phi \mu}=\sum_{i=1}^n |\inner{\phi_i}{\mu_i}|^2.
\end{equation}
 
Maximizing  $R_{\phi \mu}$ may be of interest in various
applications. For example, in a matched-filter detection problem
considered in \cite{EO01}, $R_{\phi \mu}$ represents the total output
signal-to-noise ratio. As another example, in 
a multiuser detection problem
considered in \cite{EC01}, maximizing $R_{\phi \mu}$ has the effect of
minimizing the 
multiple-access interference at the input to the proposed detector.

To obtain a more convenient expression for $R_{\phi \mu}$,
let $\Phi$ and $M$ denote  the
matrices  of columns $\phi_i$ and $\mu_i$, respectively.
Since the vectors $\mu_i$ form a normalized tight frame for $\HH$, $M$
satisfies 
\begin{equation}
\label{eq:constm}
MM^*=I_m.
\end{equation}
From Corollary~\ref{cor:gusvd}, $\Phi$ has an  SVD of the form
$\Phi=U\Sigma \F^*$, where $U$ is 
unitary, $\F$ is the FT matrix over the additive group $Q$ to which
$\Q$ is isomorphic, and
 $\Sigma$ is an $m \times n$ diagonal
matrix with diagonal elements $\sigma_i>0$.
From (\ref{eq:constm}) it follows that $M$ can be written as 
$M=U\tilde{I}Z^*$ where $Z$ is an arbitrary 
unitary matrix and $\tilde{I}$ is an $m \times n$ diagonal
matrix with diagonal elements all equal to $1$.

Let $f_i$ and $z_i$ denote the
columns of $\F^*$ and $Z^*$, respectively. Then we can express $R_{\phi \mu}$
as 
\begin{equation}
R_{\phi \mu}=\sum_{i=1}^n |\inner{\phi_i}{\mu_i}|^2=
\sum_{i=1}^n |\inner{U^*\phi_i}{U^*\mu_i}|^2=
\sum_{i=1}^n |\inner{f_i}{\overline{\Sigma}z_i}|^2,
\end{equation}
where $\overline{\Sigma}$ is an $n \times n$ diagonal matrix with the
first $m$ diagonal elements equal to $\sigma_i$, and the remaining
diagonal elements are all equal to $0$.

Our problem then reduces to finding a set of orthonormal vectors
$z_i$ that maximize $\sum_i |\inner{f_i}{\overline{\Sigma}z_i}|^2$, where
the vectors $f_i$ are also orthonormal.
Using the Cauchy-Schwarz inequality we have that
\begin{equation}
\label{eq:cs}
R_{\phi \mu} =\sum_{i=1}^n |\inner{\overline{\Sigma}^{1/2}
f_i}{\overline{\Sigma} ^{1/2}
z_i}|^2 
\leq \sum_{i=1}^n
\inner{f_i}{\overline{\Sigma}f_i}\inner{z_i}{\overline{\Sigma} z_i},  
\end{equation}
with equality if and only if $\overline{\Sigma}^{1/2} f_i=c_i
\overline{\Sigma}^{1/2} z_i$ for 
some $c_i$. In particular, we have equality for $z_i=f_i$. 
Since the components of the vectors $f_i$ all have
equal magnitude $1/\sqrt{n}$,
$\inner{f_i}{\overline{\Sigma} f_i}=1/n\sum_{k=1}^m \sigma_k \deft \alpha$
for all 
$i$, and  (\ref{eq:cs}) reduces to
\begin{equation}
R_{\phi \mu} \leq 
\alpha \sum_{i=1}^n \inner{z_i}{\overline{\Sigma}z_i} 
 =  \alpha\tr(Z^*\overline{\Sigma}Z) 
 =  \alpha \tr(\overline{\Sigma}) 
 =  n\alpha^2,
\end{equation}
with equality if $f_i=z_i$.

The normalized  tight frame vectors that maximize $R_{\phi \mu}$ are
then the columns of  
$\widehat{M}=U\tilde{I} \F^*=S^{-1/2}\Phi$ where $S=\Phi\Phi^*$, and are
equal to the canonical tight frame vectors.

\newpage
%%%%%%%%%%%%%%%%%%%%%%%%%%%%%%%%%%%%%%%%
\section{Example Of A GU Frame}
\label{sec:example}
%%%%%%%%%%%%%%%%%%%%%%%%%%%%%%%%%%%%%%
We now consider an example demonstrating the ideas of the previous
section.  

Consider the frame vectors $\phi_1=1/2[\sqrt{3}\,\,\, -1]^*,
\phi_2=1/2[\sqrt{3}\,\,\, 1]^*,
\phi_3=1/2[-\sqrt{3}\,\,\, 1]^*,
\phi_4=1/2[-\sqrt{3}\,\,\, -1]^*$, depicted in Fig.~\ref{fig:gu}.

The corresponding Gram matrix is given by
\begin{equation}
\label{eq:Gex}
G=\left[
\begin{array}{rrrr}
1 & 0.5 & -1 & -0.5 \\
0.5 & 1 & -0.5 & -1 \\
-1 & -0.5 & 1 & 0.5 \\
-0.5 & -1 & 0.5 & 1
\end{array}
\right],
\end{equation}
which is  a permuted matrix with $G=G^T$. From
Proposition~\ref{prop:permuted} it follows that the vectors $\phi_i$ are
GU. Since the vectors $\phi_i$ also span $\R^2$,
these vectors form a GU frame for $\R^2$.
\setlength{\unitlength}{.2in}
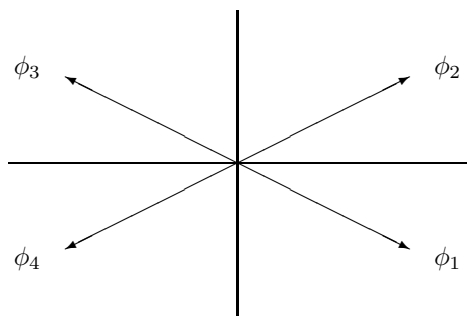
\begin{figure}[h]
\begin{center}
\begin{picture}(14,10)(0,0)
%\put(0,0){\framebox(14,10){}}

\put(0,-1){

\put(7,6){\line(0,1){4}}
\put(7,6){\line(0,-1){4}}
\put(7,6){\line(-1,0){6}}
\put(7,6){\line(1,0){6}}

\put(7,6){\vector(2,1){4.5}}
\put(12.5,8.5){\makebox(0,0){\small $\phi_2$}}
\put(7,6){\vector(2,-1){4.5}}
\put(12.5,3.5){\makebox(0,0){\small $\phi_1$}}
\put(7,6){\vector(-2,1){4.5}}
\put(1.5,8.5){\makebox(0,0){\small $\phi_3$}}
\put(7,6){\vector(-2,-1){4.5}}
\put(1.5,3.5){\makebox(0,0){\small $\phi_4$}}

}

\end{picture}
\caption{Example of a GU frame.}
\label{fig:gu}
\end{center}
\end{figure}

The vectors $\phi_i$ can be expressed as 
$\{\phi_i=U_i\phi,\,\,1 \leq i \leq
4\}$,  where $\phi=\phi_1$ and the matrices
 $\{U_i,1 \leq i
\leq 4\}$ are unitary, form an abelian group $\Q$, and are given by
\begin{equation}
U_1=I_2,\,\,\,
U_2=\left[
\begin{array}{rr}
1 & 0 \\
0 & -1 
\end{array}
\right],\,\,\,
U_3=\left[
\begin{array}{rr}
-1 & 0 \\
0 & -1 
\end{array}
\right],\,\,\,
U_4=\left[
\begin{array}{rr}
-1 & 0 \\
0 & 1 
\end{array}
\right].
\end{equation}

The multiplication table of the group $\Q$ is
\begin{equation}
\label{eq:utable}
\begin{array}{c|cccc}
 & U_1 & U_2 & U_3 & U_4 \\
\hline
U_1 & U_1 & U_2 & U_3 & U_4 \\
U_2 & U_2 & U_1 & U_4 & U_3 \\
U_3 & U_3 & U_4 & U_1 & U_2 \\
U_4 & U_4 & U_3 & U_2 & U_1.
\end{array}
\end{equation}
If we define the correspondence
\begin{equation}
\label{eq:coress}
U_1 \lra (0,0),\,\,U_2 \lra (0,1),\,\, U_3 \lra (1,0),\,\,U_4 \lra (1,1),
\end{equation}
then this table becomes the
addition table of $Q=\Z_2 \times \Z_2$:
\begin{equation}
\label{eq:atable}
\begin{array}{c|cccc}
 & (0,0) & (0,1) & (1,0) & (1,1) \\
\hline
(0,0) & (0,0) & (0,1) & (1,0) & (1,1) \\
(0,1) & (0,1) & (0,0) & (1,1) & (1,0) \\
(1,0) & (1,0) & (1,1) & (0,0) & (0,1) \\
(1,1) & (1,1) & (1,0) & (0,1) & (0,0).
\end{array}
\end{equation}
Only the way in which the elements are labeled distinguishes the table of
(\ref{eq:atable}) from the table of (\ref{eq:utable});
thus $\Q$ is isomorphic to $Q$.
Over $Q=\Z_2 \times \Z_2$, the FT matrix $\F$ is the
Hadamard matrix
\begin{equation}
\F=\frac{1}{2}\left[
\begin{array}{rrrr}
1 & 1 & 1 & 1 \\
1 & -1 & 1 & -1 \\
1 & 1 & -1 & -1 \\
1 & -1 & -1 & 1
\end{array}
\right].
\end{equation}

From Theorem~\ref{thm:gu} the dual frame vectors and
the canonical tight frame vectors are also GU with generators $\bar{\phi}$
and $\mu$ respectively, whose equations are given in the theorem.
Thus, to compute the dual and canonical tight frame vectors we compute these
generators and then apply the group $\Q$.

We first determine
the FT of the first row of $G$ denoted by $s$:
\begin{equation}
\hat{s}=\F s=\frac{1}{2}[0\,\,0\,\,3\,\,1]^*. \label{shat}
\end{equation} 
Using Theorem~\ref{thm:gu}, it follows from (\ref{shat}) that the frame bounds 
are given by $A=1$ and $B=3$.
Next, we compute the vectors $\hat{\phi}(h)$ which are the columns of 
\begin{equation}
\Phi\F=\left[
\begin{array}{rrrr}
0 & 0 & 1.7 & 0 \\
0 & 0 & 0 & -1 \\
\end{array}
\right].
\end{equation} 
Using the equations of the theorem, we then have that
$\bar{\phi}=[0.3\,\, -0.5]^*$ and $\mu=[0.5\,\, -0.5]^*$.
By applying the group $\Q$ to these generators we obtain that the dual
frame vectors are the columns of 
\begin{equation}
\label{eq:oPhi}
\overline{\Phi}=\left[
\begin{array}{rrrr}
0.3 & 0.3 & -0.3 & -0.3 \\
-0.5 & 0.5 & 0.5 & -0.5 \\
\end{array}
\right],
\end{equation} 
and the canonical tight
frame vectors are the columns of 
\begin{equation}
\label{eq:M}
M=\left[
\begin{array}{rrrr}
0.5 & 0.5 & -0.5 & -0.5 \\
-0.5 & 0.5 & 0.5 & -0.5 \\
\end{array}
\right].
\end{equation}

Comparing (\ref{eq:oPhi}) and (\ref{eq:M}) with the original frame
vectors $\phi_i$, it is
evident that the dual and canonical tight frame vectors have the same
symmetries as the original frame vectors, as illustrated in 
Fig.~\ref{fig:dc}.

\setlength{\unitlength}{.2in}
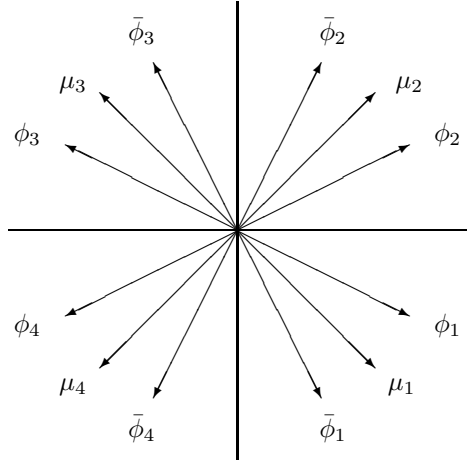
\begin{figure}[h]
\begin{center}
\begin{picture}(14,14)(0,0)
%\put(0,0){\framebox(14,14){}}

\put(0,1){

\put(7,6){\line(0,1){6}}
\put(7,6){\line(0,-1){6}}
\put(7,6){\line(-1,0){6}}
\put(7,6){\line(1,0){6}}

\put(7,6){\vector(2,1){4.5}}
\put(12.5,8.5){\makebox(0,0){\small $\phi_2$}}
\put(7,6){\vector(2,-1){4.5}}
\put(12.5,3.5){\makebox(0,0){\small $\phi_1$}}
\put(7,6){\vector(-2,1){4.5}}
\put(1.5,8.5){\makebox(0,0){\small $\phi_3$}}
\put(7,6){\vector(-2,-1){4.5}}
\put(1.5,3.5){\makebox(0,0){\small $\phi_4$}}

\put(7,6){\vector(1,2){2.2}}
\put(9.5,11.2){\makebox(0,0){\small $\bar{\phi}_2$}}
\put(7,6){\vector(1,-2){2.2}}
\put(9.5,.8){\makebox(0,0){\small $\bar{\phi}_1$}}
\put(7,6){\vector(-1,2){2.2}}
\put(4.5,11.2){\makebox(0,0){\small $\bar{\phi}_3$}}
\put(7,6){\vector(-1,-2){2.2}}
\put(4.5,.8){\makebox(0,0){\small $\bar{\phi}_4$}}

\put(7,6){\vector(1,1){3.6}}
\put(11.5,9.8){\makebox(0,0){\small $\mu_2$}}
\put(7,6){\vector(1,-1){3.6}}
\put(11.3,1.9){\makebox(0,0){\small $\mu_1$}}
\put(7,6){\vector(-1,1){3.6}}
\put(2.7,9.8){\makebox(0,0){\small $\mu_3$}}
\put(7,6){\vector(-1,-1){3.6}}
\put(2.7,1.9){\makebox(0,0){\small $\mu_4$}}

}

\end{picture}
\caption{Symmetry property of the frame vectors
${\phi_i}$, the dual frame vectors
vectors $\bar{\phi_i}$, and the canonical tight frame
vectors ${\mu_i}$.
$\bar{\phi}_i$ are the columns of $\overline{\Phi}$ given by
(\ref{eq:oPhi}), and 
$\mu_i$ are the columns of $M$ given by (\ref{eq:M}).
The frame vectors, dual frame vectors, and the canonical tight frame vectors
all have the same symmetry properties.}
\label{fig:dc}
\end{center}
\end{figure}

\newpage
%%%%%%%%%%%%%%%%%%%%%%%%%%%%%%%%%%%%%%%%
\section{Compound GU Frames}
\label{sec:cguf}
%%%%%%%%%%%%%%%%%%%%%%%%%%%%%%%%%%%%%%

In Section~\ref{sec:dual} we showed that the dual and canonical tight frame
vectors associated with a GU frame are themselves GU and can therefore
be computed using a single
generator. 
In this section,
we consider a class of frames which consist of subsets that are GU,
and are therefore referred to as {\em compound geometrically uniform
(CGU) frames}. As we show, the dual and canonical tight frame vectors associated
with a CGU frame share the same symmetries as the original frame
and can be computed using a {\em set} of generators.

A set of frame vectors 
$\{\phi_{ik},1 \leq i \leq l,1 \leq k\leq r\}$ is CGU if $\phi_{ik}=U_i\phi_k$ 
for  some generating vectors
$\{\phi_k,1 \leq k \leq r \}$, and the matrices
$\{U_i,1 \leq i \leq l\}$ are unitary and form an abelian group $\Q$.

A CGU frame is in general not GU. However,
for every $k$, the vectors $\{\phi_{ik},1 \leq i \leq l\}$ are a GU
vector set with generating group $\Q$.

A special case of CGU frames are  filter bank frames studied in
\cite{B97,BHF98,Cve98}, in which $\Q$ is the group of translations 
by  integer multiples of the subsampling factor, and the generating
vectors are  the filter bank synthesis filters.  We may therefore view
CGU frames as an extension of filter bank frames to the more general group 
case.

As we show in the following proposition, the frame bounds of a CGU
frame can be bounded by
the sum of the norms of the generating vectors. 
\begin{prop}
\label{prop:boundsg}
Let $\SSS = \{\phi_{ik} = U_i\phi_k, 1 \leq i \leq l,1 \leq k \leq r\}$ be
a compound 
geometrically 
uniform frame with frame bounds $A$ and $B$, where $\{\phi_k,1 \leq k
\leq r\}$ is an
arbitrary set of generating vectors. Then
$A\leq \frac{l}{m}\sum_{k=1}^r||\phi_k||^2 \leq B$. 
If in addition the frame is tight, then
$A=\frac{l}{m} \sum_{k=1}^r||\phi_k||^2$. 
\end{prop}
\begin{proof}
We can express the frame operator corresponding to the frame vectors
$\phi_i$ as
\begin{equation}
S=\sum_{i=1}^l\sum_{k=1}^r U_i\phi_k\phi_k^*U_i^*.
\end{equation}
Then,
\begin{equation}
\tr(S)=\sum_{i=1}^l\sum_{k=1}^r \tr(U_i\phi_k\phi_k^*U_i^*)=
l\sum_{k=1}^r ||\phi_k||^2,
\end{equation}
so that
\begin{equation}
\sum_{i=1}^m \lambda_i(S)=\tr(S)=l\sum_{k=1}^r||\phi_k||^2.
\end{equation}
Therefore,
\begin{equation}
A=\min_i \lambda_i(S) \leq \frac{1}{m}\sum_{i=1}^m
\lambda_i(S)=\frac{l}{m}\sum_{k=1}^r ||\phi_k||^2,
\end{equation}
and
\begin{equation}
B=\max_i \lambda_i(S) \geq \frac{1}{m}\sum_{i=1}^m
\lambda_i(S)=\frac{l}{m}\sum_{k=1}^r ||\phi_k||^2.
\end{equation}
\end{proof}

%%%%%%%%%%%%%%%%%%%%%%%%%%%%%%%%%%%%%%%%%%%%%%%%%
\subsection{Dual and canonical tight frames associated with CGU frames}
\label{sec:dualcgu}
%%%%%%%%%%%%%%%%%%%%%%%%%%%%%%%%%%%%%%%%%%%%%%%%

We now show that the dual and canonical tight frames associated with a CGU
frame are also CGU.

Expressing the frame operator as
\begin{equation}
S=\sum_{i=1}^l\sum_{k=1}^r \phi_{ik}\phi_{ik}^*=
\sum_{i=1}^lU_i \bl \sum_{k=1}^r \phi_{k}\phi_{k}^* \br U_i^*,
\end{equation}
for all $j$ we have that,
\begin{eqnarray}
SU_j & = & \sum_{i=1}^lU_i \bl \sum_{k=1}^r \phi_{k}\phi_{k}^* \br U_i^*U_j
\nonumber \\
& = & U_j\sum_{i=1}^lU_j^*U_i \bl \sum_{k=1}^r \phi_{k}\phi_{k}^* \br U_i^*U_j
\nonumber \\
& = & U_j  \sum_{i=1}^lU_i \bl \sum_{k=1}^r \phi_{k}\phi_{k}^* \br U_i^*
\nonumber \\
& = & U_jS,
\end{eqnarray}
since $\{U_j ^*U_i,1 \leq i \leq l\}$ is just a permutation of $\Q$.
Thus $S$ commutes with $U_j$, so that  $S^{-1}$ and $S^{-1/2}$ also commute
with $U_j$ for all $j$. Then, the dual frame vectors $\bar{\phi}_{ik}$
of the vectors $\phi_{ik}$ are given by
\begin{equation}
\bar{\phi}_{ik}=S^{-1}\phi_{ik}=S^{-1}U_i\phi_k=U_iS^{-1}\phi_k=U_i 
\bar{\phi}_k,
\end{equation}
where $\bar{\phi}_k=S^{-1} \phi_k$, which shows that the dual frame vectors
$\{\bar{\phi}_{ik}=S^{-1}\phi_{ik}\}$ are CGU with generating group equal
to $\Q$.

Similarly,
\begin{equation}
\mu_{ik}=S^{-1/2}\phi_{ik}=S^{-1/2}U_i\phi_k=U_iS^{-1/2}\phi_k=U_i \mu_k,
\end{equation}
where $\mu_k=S^{-1/2} \phi_k$, which shows that the canonical tight frame vectors 
$\{\mu_{ik}=S^{-1/2}\phi_{ik}\}$ are also CGU with generating
group $\Q$.

Therefore, to compute the dual frame vectors or the canonical tight frame
vectors all we
need is to compute  the generating vectors $\{\bar{\phi}_k,1 \leq k
\leq r\}$ and $\{\mu_k,1 \leq k \leq r\}$,
respectively. 
The remaining frame vectors are then obtained by applying the group $\Q$ to the
corresponding set of generating vectors. 

%Since the frame is not generally GU, we can no longer
%compute these generating vectors efficiently using the FT matrix. 
%In practice, to avoid computing $S^{-1}$ and $S^{-1/2}$ we may choose to
%approximate the generating vectors using the series expansions
%(\ref{eq:sinv}) and (\ref{eq:ssqrt}).

%For example, to compute $S^{-1}$ approximately we can truncate  
%the series expansion (\ref{eq:sinv}).
%Estimates of the resulting reconstruction error are given in \cite{D92}.
%Suppose we use a first order approximation
%obtained by retaining only the 
%$i=0$ term,
%\begin{equation}
%(S^{-1})^{(0)}=\frac{2}{A+B}I.
%\end{equation}
%Then the  approximate dual frame generators are given by
%\begin{equation}
%\label{eq:app}
%\bar{\phi}_{k} \approx\frac{2}{A+B}\phi_{k}.
%\end{equation}
%If we now 
%reconstruct $x \in \HH$ from the frame coefficients
%$\inner{\phi_{ik}}{x}$ using 
%the approximate dual frame vectors resulting from (\ref{eq:app}),
%then the reconstruction error incurred by 
%this approximation can be bounded 
%in terms of the frame 
%bounds $A$ and $B$. With $\hat{x}$ denoting the reconstructed vector
%we have \cite{D90}
%\begin{equation}
%\|\hat{x}-{x}\| \leq \frac{B/A-1}{B/A+1} \| x\|.
%\end{equation}
%This reconstruction error is small if the tightest possible frame
%bounds $A$ and $B$ 
%satisfy $B/A\,\approx\,1$, in which case  the underlying frame is
%called a snug frame.

%Starting from (\ref{eq:ssqrt}), we  can approximate the canonical dual frame
%vectors in a similar manner.

%%%%%%%%%%%%%%%%%%%%%%%%%%%%%%%%%%%%%%%%
\subsection{CGU frames with GU generators}
%%%%%%%%%%%%%%%%%%%%%%%%%%%%%%%%%%%%%%

A special class of CGU frames  is {\em CGU
frames with GU generators} in which the generating vectors $\{\phi_k,1
\leq k \leq r\}$ are
themselves GU. Specifically, $\{\phi_k=V_k \phi\}$ for some generator $\phi$,
where the matrices $\{V_k,1 \leq k \leq r\}$ are unitary, and form an
abelian group $\G$.

Now suppose that $U_p$ and $V_t$ commute up to a phase factor for all
$t$ and $p$ so that $U_pV_t=V_tU_pe^{j\theta(p,t)}$ where $\theta(p,t)$
is an arbitrary phase function that may depend on the indices $p$ and
$t$. In this case we say that $\Q$ and $\G$ commute up to a phase
factor. Then for  all $p,t$,
\begin{eqnarray}
SU_pV_t & = & \sum_{i=1}^lU_i \bl \sum_{k=1}^rV_k\phi \phi^*V_k^* \br 
U_i^*U_pV_t
\nonumber \\
& = & U_pV_t\sum_{i=1}^lV_t^*U_p^*U_i \bl \sum_{k=1}^rV_k \phi \phi^*V_k^*\br
U_i^*U_pV_t 
\nonumber \\
& = & U_p V_t \sum_{i=1}^lV_t^*U_i^* \bl \sum_{k=1}^rV_k \phi \phi^*V_k^*\br
U_i V_t 
\nonumber \\
& = & U_p V_t \sum_{i=1}^lU_i^* \bl \sum_{k=1}^rV_t^* V_k \phi 
\phi^*V_k^*V_t\br U_i
\nonumber \\
& = & U_p V_t \sum_{i=1}^lU_i \bl \sum_{k=1}^rV_k\phi \phi^*V_k^*\br U_i^* 
\nonumber \\
& = & U_pV_t S.
\end{eqnarray}
The dual frame vectors $\bar{\phi}_{ik}$
of the vectors $\phi_{ik}$ are then given by
\begin{equation}
\bar{\phi}_{ik}=S^{-1}\phi_{ik}=S^{-1}U_iV_k \phi=U_iV_kS^{-1} \phi=U_i 
V_k\bar{\phi},
\end{equation}
where $\bar{\phi}=S^{-1}\phi$.
Similarly,
\begin{equation}
\mu_{ik}=S^{-1/2}\phi_{ik}=S^{-1/2}U_iV_k\phi=U_iV_kS^{-1/2}\phi=U_iV_k
\mu
\end{equation}
where $\mu=S^{-1/2}\phi$. Thus even though the frame is not in
general GU, the dual and canonical tight frame vectors can be computed using
a single generating vector. 

As we now show, in the special case in which $\theta=0$ so that
$U_iV_k=V_kU_i$ for 
all $i,k$, the resulting frame is GU. To this end we need to show
that the unitary matrices $\Q'=\{Q_{ik}=U_iV_k\}$ form an abelian group. 
First, $Q_{ik}Q_{jt}=U_iV_kU_jV_t= U_iU_jV_kV_t$. Since $U_iU_j \in
\Q$ and $V_kV_t \in \G$, 
$Q_{ik}Q_{jt} \in \Q'$. 
Next, $Q_{ik}^{-1}=V_k^*U_i^*=U_i^*V_k^* \in \Q'$ since  $U_i^* \in
\Q$ and $V_k^* \in \G$. Also, $I \in \Q'$ since $I \in \Q$ and
$I \in \G$. Finally,
$Q_{ik}Q_{jt}=U_iV_kU_jV_t=U_iU_jV_tV_k=U_jV_tU_iV_k=Q_{jt}Q_{ik}$. 

A special case of CGU frames with GU generators for which  $\Q$ and $\G$
commute up to a phase factor 
are Weyl-Heisenberg (WH) frames
\cite{D92,HW89,fei_book}. If the WH frame is critically sampled, then 
$\theta(p,t)=0$ and the WH frame reduces to a GU frame. In the more
general oversampled case,
$\theta(p,t) \neq 0$. CGU frames with GU generators can
hence be viewed as a generalization of WH frames to the
group case.
 
To summarize, we have the following theorem:
\begin{theorem}[Compound GU frames]
\label{thm:cgu}
Let $\SSS = \{\phi_{ik} = U_i\phi_k, 1 \leq i \leq l,1 \leq k \leq r\}$ be
a compound 
geometrically 
uniform frame with  generating vectors $\{\phi_k,1 \leq k \leq r\}$
and generating group $\Q$,
and let $S$ be the frame operator corresponding to the frame
vectors $\{\phi_{ik}\}$.
Then
\begin{enumerate}
\item the dual frame vectors $\{\bar{\phi}_{ik},1 \leq i \leq l,1 \leq 
k \leq r\}$ are
compound geometrically uniform with generating group $\Q$ and
generating vectors $\{\bar{\phi}_k=S^{-1}\phi_k,1 \leq k \leq r\}$;
\item the canonical tight frame vectors $\{\mu_{ik},1 \leq i \leq l,1 \leq 
k \leq r\}$ are
compound geometrically uniform with generating group $\Q$ and
generating vectors $\{\mu_k=S^{-1/2}\phi_k,1 \leq k \leq r\}$.
\end{enumerate}
If in addition the generating vectors $\{\phi_k=V_k\phi,1 \leq k
\leq r\}$ are geometrically uniform
with  $U_iV_k=V_kU_ie^{j\theta(i,k)}$
for all 
$i,k$, then
\begin{enumerate}
\item $\bar{\phi}_{ik}=U_iV_k\bar{\phi}$ where $\bar{\phi}=S^{-1}\phi$;
\item  $\mu_{ik}=U_iV_k\mu$ where $\mu=S^{-1/2}\phi$;
\item if in addition $\theta(i,k)=0$ for all $i,k$, then the vectors
$\{\phi_{ik}, 1 \leq i \leq l,1 \leq k \leq r\}$ form a geometrically
uniform frame. 
\end{enumerate}
\end{theorem}

\newpage 
%%%%%%%%%%%%%%%%%%%%%%%%%%%%%%%%%%%%%
\section{Pruning GU Frames}
\label{sec:pruning}
%%%%%%%%%%%%%%%%%%%%%%%%%%%%%%%%%%%%%%

In applications it is often desirable to know how a frame behaves when one
or more frame elements are removed. In particular, it is important to
know or to be able to estimate the frame bounds of the reduced frame.
In general, if no structural constraints are imposed on a frame
this behavior will depend critically on the particular frame elements
removed. For example, removal of a particular frame element may
destroy the frame 
property so that the remaining vectors do not constitute a frame anymore, 
whereas if a different  element is removed
the remaining vectors may still constitute a frame.

One of the prime applications of frames is signal analysis and
synthesis, where a  signal is expanded by computing the inner products of the signal with
the frame elements.
The resulting coefficients  are subsequently stored,
transmitted, quantized or manipulated in some way. 
In particular,  a coefficient may be lost (\eg due to a transmission error)
 which results in a reconstructed signal that is
equivalent to an expansion using  a pruned
frame obtained by 
removing the corresponding frame vector. 

Recently, there has been increased interest in using frames for
multiple-description source coding where a signal is expanded into a
redundant set of functions and the resulting coefficients are
transmitted over a lossy packet network where one or more of the
coefficients can be lost because a packet is dropped \cite{GKV98}. The
goal of multiple description source coding is to ensure a gradually
behaving reconstruction quality as a function of the number of dropped
packets. When using  frames in this context, the reconstruction quality is
often governed by the frame bound ratio of the pruned frame. 
If the packets are dropped with equal probability, then it is desirable
that  the frame
bound ratio should deteriorate uniformly irrespectively of the
particular frame element that is removed.  In the following, we show that
GU frames have this property. We will furthermore demonstrate that if
the original frame is  a tight GU frame, then the frame bound ratio of
the pruned 
frame obtained by removing one frame element can be computed
exactly. We  also consider the case where sets of frame elements
are removed.

\begin{theorem}[Pruned GU frames]
\label{thm:pruned_gu}
Let $\SSS = \{\phi_i = U_i\phi, U_i \in \Q\}$ be a geometrically
uniform frame generated by a finite abelian group $\Q$ of unitary
matrices, where $\phi$ is an arbitrary generating vector. Let
$\Phi$ be the matrix 
of
columns $\phi_i$, and let $S=\Phi\Phi^{\ast}$ be the corresponding frame
operator.  
Let $\SSS(j) = \{\phi_i = U_i\phi, U_i \in \Q, i\,\neq\,j\}$
be the pruned set obtained by removing the element $\phi_{j}$. 
Then the eigenvalues of the frame operator corresponding to the pruned set
do not depend on the particular element $\phi_{j}$ removed.
\end{theorem}
\begin{proof}
The frame operator corresponding to the pruned frame is given by
\begin{equation}
S(j)=\sum_{i=1}^{n}U_i\phi \phi^{\ast} U^{\ast}_i-U_j\phi \phi^{\ast} 
U^{\ast}_{j}.
\end{equation}
Since $U_j$ is unitary, the eigenvalues of $S(j)$ are equal to the
eigenvalues of 
$U_j^*S(j)U_j$. But,
\begin{equation}
U^{\ast}_{j} S(j) U_{j} = U^{\ast}_{j}\sum_{i=1}^{n}U_i\phi \phi^{\ast} 
U^{\ast}_i U_{j}
-\phi \phi^{\ast}
=\sum_{i=1}^{n}U_i\phi \phi^{\ast} U^{\ast}_i - \phi \phi^{\ast} = S - \phi \phi^{\ast}.
\end{equation}
Since $U_j^*S(j)U_j$ is independent of $j$, the eigenvalues of $S(j)$
do not depend on $j$.
\end{proof}

In  general it is difficult to provide estimates on the frame bounds
of the pruned frame. 
However, in the special case where the original GU frame is a tight
frame, these bounds can be determined exactly.

\begin{cor}[Pruned Tight GU frames]
\label{cor:tight_p_gu}
Let $\SSS = \{\phi_i = U_i\phi, U_i \in \Q\}$ be a geometrically
uniform tight frame generated by a finite abelian group $\Q$ of unitary
matrices, where $\phi$ is a unit norm generating vector. Let
$\Phi$ be the matrix 
of
columns $\phi_i$, and let $S=\Phi\Phi^{\ast}$ be the corresponding
frame operator.  
Let  $\SSS(j) = \{\phi_i = U_i\phi, U_i \in \Q, i\,\neq\,j\}$
be the pruned set obtained by removing the element $\phi_{j}$. 
Then the eigenvalues of the frame operator corresponding to the pruned set
are given by $\lambda_{1}=\frac{n}{m}-1$ and
$\lambda_{i}=\frac{n}{m},\,2 \leq i \leq n$,
independent of  $\phi_{j}$. 
\end{cor}
\begin{proof}
Since $\SSS$ is a tight frame with $||\phi||=1$, from
Proposition~\ref{prop:bounds} the frame bound $A=n/m$ and $S=(n/m)I_m$.
Then,
\begin{equation}
S(j)=\frac{n}{m}I_m-U_j\phi \phi^{\ast} U^{\ast}_{j},
\end{equation}
and
\begin{equation}
\label{eq:usju}
U^{\ast}_{j} S(j) U_{j} =  \frac{n}{m}I_m-\phi \phi^{\ast}.
\end{equation}
Since $||\phi||=1$, $\phi\phi^*$ has one eigenvalue equal to $1$, and
the remaining eigenvalues equal to $0$. 
The eigenvalues of $S(j)$ are therefore given by
$\lambda_{1}=\frac{n}{m}-1$ and
$\lambda_{i}=\frac{n}{m},\,2 \leq i \leq n$.
\end{proof}

An immediate consequence of Corollary~\ref{cor:tight_p_gu} is that the
frame bound ratio of the pruned frame is given by $B/A=1/(1-m/n)$,
which is close to $1$ for large redundancy $r=n/m$.

%We conclude that GU frames are particularly attractive in multiple
%description applications since the frame bound
%ratio of the pruned frame is independent of the particular frame
%element removed, and for large $n/m$ this bound 
%will be close to $1$ resulting in
%good numerical properties of the reconstruction.

We next consider the case where multiple frame elements are removed. 
\begin{cor}
\label{cor:mult_tight_pruned_gu}
Let $\SSS = \{\phi_i = U_i\phi, U_i \in \Q\}$ be a geometrically
uniform frame generated by a finite abelian group $\Q$ of unitary
matrices, let $\Phi$ be the matrix of columns $\phi_i$, and let
$S=\Phi\Phi^{\ast}$ be the corresponding frame operator.  Let ${\cal J}$ be
a  set of indices, and let ${\cal J}(k)$
denote the  set of indices $i$ such that $U_i= U_k U_j$ for fixed 
$k$ and 
$j\,\in\,{\cal J}$.  
Let $\SSS(k) = \{\phi_i = U_i\phi,
U_i \in \Q, i\,\neq\,{\cal J}(k)\}$ be a pruned set obtained by removing
the elements $\phi_{i}$ with $i\,\in\,{\cal J}(k)$. 
Then the eigenvalues of the
frame operator corresponding to the pruned set are independent of $k$.
\end{cor}
\begin{proof}
The  frame operator corresponding to the pruned frame is
given by
\begin{equation}
S(k)=\sum_{i=1}^{n}U_i\phi \phi^{\ast} U^{\ast}_i-U_k \left[\sum_{j\,\in\,{\cal 
J}}U_j\phi \phi^{\ast} U^{\ast}_{j} \right]U^{\ast}_k.
\end{equation}
Then,
\begin{eqnarray*}
U^{\ast}_{k} S(k) U_{k} & = & U^{\ast}_{k}\sum_{i=1}^{n}U_i\phi \phi^{\ast} 
U^{\ast}_i U_{k}
-\sum_{j\,\in\,{\cal J}}U_j\phi \phi^{\ast} U^{\ast}_{j}\\
&& = \sum_{i=1}^{n}U_i\phi \phi^{\ast} U^{\ast}_i - \sum_{j\,\in\,{\cal 
J}}U_j\phi \phi^{\ast} U^{\ast}_{j}
\end{eqnarray*}
is independent of $k$, and consequently the eigenvalues of $S(k)$ do not
depend on $k$.
\end{proof}

To conclude this section, GU frames have  strong symmetry properties
in the sense that removing any one of the elements leads to a vector
set with  bounds independent
of the particular element removed. Moreover, if the original frame is tight, 
then we can compute the
bounds of the pruned frame exactly. 

%%%%%%%%%%%%%%%%%%%%%%%%%%%%%%%%%%%%%%%%
\section{Constructing GU Frames}
\label{sec:optimal}
%%%%%%%%%%%%%%%%%%%%%%%%%%%%%%%%%%%%%%

Suppose we are given a
 set of vectors $\{\varphi_i,1
\leq i \leq n\}$, that form a frame for an $m$-dimensional space
$\HH$. We would like to construct a geometrically uniform frame
$\{\phi_i,1 \leq i \leq n\}$ from
the vectors $\{\varphi_i\}$.

From
Theorem~\ref{thm:guft} it follows that the vectors $\{\phi_i\}$ form
a  GU frame if and only if the Gram matrix $G$ has rank $m$, and 
is diagonalized by a FT matrix $\F$ over a finite product of
cyclic groups. 
There are many ways to construct a frame from a given set of
frame vectors $\varphi_i$ that satisfy these properties. For example,
let $F$ be the matrix of columns $\varphi_i$, and let $F$ have an SVD
$F=Q\Lambda V^*$, where $\Lambda$ is a diagonal matrix with diagonal
elements $\lambda_i$.  Then  the columns of 
$\Phi=F V \Sigma \F^*=Q\Lambda\Sigma \F^*$ form a GU frame, where $\F$
is any FT matrix over a 
product of cyclic groups, and  $\Sigma$ is an arbitrary
diagonal matrix with 
diagonal elements $\sigma_i>0$.
The frame bounds of the resulting GU frame are given
by $A=\min_i \lambda_i^2\sigma_i^2$ and $B=\max_i
\lambda_i^2\sigma_i^2$, so that we can choose the  diagonal
matrix $\Sigma$ to control these  bounds.
In particular, choosing $\Sigma=I$ we have that the
columns of $\Phi=FV\F^*=Q\Lambda\F^*$ form a GU frame. This choice has the
property that the frame bounds of the GU frame are equal
to the frame bounds of the original frame.

We now consider the problem
of constructing an {\em optimal} GU frame.
Specifically, let $\{\varphi_i,1 \leq i \leq n\}$ denote a frame for
$\HH$, and suppose we wish to construct a GU frame $\{\phi_i\}$ from
the vectors 
$\{\varphi_i\}$. 
A reasonable approach is to find a set of vectors $\phi_i$ that span
$\HH$, and are ``closest'' to the vectors $\varphi_i$ in the least-squares
sense.  Thus we seek vectors $\phi_i$ that minimize the least-squares error
$E$, defined by
\begin{equation}
\label{eq:serror}
E=\sum_{i=1}^n \inner{e_i}{e_i},
\end{equation} 
where $e_i$ denotes the $i$th error vector
\begin{equation}
\label{eq:error}
e_i=\varphi_i-\phi_i,
\end{equation} 
subject to the constraint that the  vectors $\phi_i$ form a GU frame.

If the vectors $\phi_i$ are GU, then their Gram matrix
$G=\Phi^*\Phi$ is a permuted matrix with rank $m$, diagonalized by a
FT matrix $\F$. Thus, the inner
products $\{\inner{\phi_i}{\phi_j}\}$ must satisfy
\begin{equation}
\label{eq:guconst}
\{\inner{\phi_i}{\phi_j},1 \leq j \leq n\}=\beta^2 P_i\{a_j,1 \leq j \leq n\},
\end{equation}
where $P_i\{a_j,1 \leq j \leq n\}$ is a permutation of the
numbers $\{a_j,1 \leq j \leq n\}$,  $\beta>0$ is a scaling
factor, and the numbers $\{a_j,1 \leq j \leq n\}$ are chosen such
that the matrix $R$ whose $i$th row is equal to $P_i\{a_j,1 \leq j
\leq n\}$ is Hermitian, non-negative definite and diagonalized by a FT
matrix $\F$.

In our development of the optimal GU frame vectors we assume that the
permutations $P_i$ in (\ref{eq:guconst}) are specified.
Since these  permutations determine the additive  group
$Q$  over which  the
FT matrix $\F$ is defined, we assume that $\F$ is specified.
We then consider three different constraints on the vectors $\phi_i$. 
First we consider the case in which both the numbers
$\{a_j,1 \leq j \leq n\}$ and the scaling $\beta$ in (\ref{eq:guconst})
are known. The  GU frame 
minimizing  the least-squares error $E$ of (\ref{eq:serror})-(\ref{eq:error})
subject to this constraint is derived in Section~\ref{sec:scgu}, and 
is referred to as the scaled-constrained
least-squares GU frame (SC-LSGUF).  
Next, we  consider
the case in which the numbers $\{a_j,1 \leq j \leq n\}$ in
(\ref{eq:guconst}) are known, and the scaling $\beta$ is chosen to
minimize  $E$. The
resulting GU frame is referred to as  
the constrained least-squares GU frame (C-LSGUF), and is derived in
Section~\ref{sec:cgu}. Finally,
in Section~\ref{sec:gu} we consider
the more general case in which both the numbers $\{a_j,1 \leq j
\leq n\}$ and the scaling $\beta$ in
(\ref{eq:guconst}) are chosen to minimize $E$.
The resulting GU frame is referred to as 
the  least-squares GU frame (LSGUF).

%%%%%%%%%%%%%%%%%%%%%%%%%%%%%%%%%%%%%%%%
\subsection{Scaled-constrained least-squares GU frame}
\label{sec:scgu}
%%%%%%%%%%%%%%%%%%%%%%%%%%%%%%%%%%%%%%

We first consider the case in which the  
 numbers $\{a_j,1 \leq j \leq n\}$ and the scaling $\beta$ in
 (\ref{eq:guconst}) are known.
Thus, we seek the set of  vectors $\{\phi_i\}$ that minimize the
least-squares error $E$ of (\ref{eq:serror})-(\ref{eq:error}) subject to
 the constraint  
\begin{equation}
\label{eq:constb}
\Phi^*\Phi=\beta_0^2 R=\beta_0^2 \F A \F^*,
\end{equation}
where $\Phi$ is the matrix of columns $\phi_i$,
$\beta_0$ is a known scaling factor, and 
$R$ is the matrix whose $i$th row is equal to 
$P_i\{a_j,1 \leq j \leq n\}$ where the numbers $\{a_j,1 \leq j \leq
n\}$ are given such that $R$ is diagonalized by $\F$, and $A$ 
is a diagonal
matrix with diagonal elements $\{\alpha_j=n^{1/2}\hat{a}_j,1 \leq j
\leq n\}$ where  
$\{\hat{a}_j,1 \leq j \leq n\}$
is the FT  of the
sequence $\{a_j,1 \leq j \leq n\}$.
From (\ref{eq:constb}), the frame bounds of the vectors $\{\phi_i\}$
are given by  
$A=\beta_0^2 \min_j \alpha_j$ and $B=\beta_0^2 \max_j \alpha_j$.

This problem has been solved in the context of general
least-squares inner product shaping \cite{E01}. The optimal SC-LSGUF
vectors $\hat{\phi}_i$ are the columns of $\widehat{\Phi}$,
given by 
\begin{equation}
\widehat{\Phi}=\beta_0 UV^* \Sigma \F^*,
\end{equation}
where $U$ and $V$ are the right-hand unitary matrix and left-hand
unitary matrix respectively in the SVD of $F\F \Sigma^*$,
$F$ is the matrix of columns $\varphi_i$, and $\Sigma$ is an $m \times
n$ diagonal matrix with diagonal elements $\sqrt{\alpha_i}$ for values
of $i$ for which $\alpha_i \neq 0$.

If $F RF^*=\beta^{2}_0 F \F A \F^* F$ is invertible, then we may express
$\widehat{\Phi}$ as  
\begin{equation}
\widehat{\Phi}=\beta_0 (F\F A\F^* F^*)^{-1/2}F\F A\F^*=\beta_0
(FRF^*)^{-1/2}FR.
\end{equation}

%%%%%%%%%%%%%%%%%%%%%%%%%%%%%%%%%%%%%%%%
\subsection{Constrained least-squares GU frame}
\label{sec:cgu}
%%%%%%%%%%%%%%%%%%%%%%%%%%%%%%%%%%%%%%

We now consider the case in which  the numbers $\{a_j,1 \leq j \leq
n\}$ are known, but the scaling $\beta$ is not specified. Thus, we
seek a set of  vectors 
$\{\phi_i\}$ that minimize the least-squares error $E$ 
subject to 
\begin{equation}
\label{eq:constbu}
\Phi^*\Phi=\beta^2 R,
\end{equation}
where $\beta>0$.
 
The least-squares error $E$ of
(\ref{eq:serror})-(\ref{eq:error}) may  be expressed as
\begin{eqnarray}
\label{eq:errorm}
E & = & \tr\bl (\Phi-F)^*(\Phi-F) \br \nonumber  \\
& = & \tr(\Phi^*\Phi)+\tr(F^*F)-2\re{\tr(F^*\Phi)} \nonumber \\
& = & \beta^2 \tr(R)+\tr(F^*F)-2\re{\tr(F^*\Phi)}.
\end{eqnarray} 
Let $\widetilde{\Phi}=(1/\beta)\Phi$. Then 
minimizing $E$ is equivalent to minimizing 
\begin{equation}
\label{eq:ep}
E'=\beta^2 \tr(R)-2\beta \re{\tr(F^*\widetilde{\Phi})},
\end{equation} 
subject to 
\begin{equation}
\label{eq:constbt}
\widetilde{\Phi}^*\widetilde{\Phi}=R.
\end{equation}

To determine the optimal matrix $\Phi$ we have to minimize $E'$ with
respect to $\beta$ and $\widetilde{\Phi}$. Fixing $\beta$ and
minimizing  with respect to $\widetilde{\Phi}$, 
the optimal value of
$\widetilde{\Phi}$ is given by the SC-LSGUF of Section~\ref{sec:scgu}
with scaling $\beta_0=1$, 
so that
\begin{equation}
\widehat{\widetilde{\Phi}}=UV^* \Sigma \F^*.
\end{equation}
If $F RF^*$ is invertible, then 
\begin{equation}
\widehat{\widetilde{\Phi}}=(F\F A\F^* F^*)^{-1/2}F\F A\F^*=
(FRF^*)^{-1/2}FR.
\end{equation}

Substituting $\widehat{\widetilde{\Phi}}$ back into (\ref{eq:ep}), and
minimizing with respect to $\beta$, the optimal value of $\beta$ is
given by
\begin{equation}
\label{eq:hatb2}
\widehat{\beta}=\frac{\re{\tr(F^*\widehat{\widetilde{\Phi}})}}{\tr(R)}=
\frac{\tr(UV^* \Sigma \F^*)}{\tr(R)},
\end{equation}  
which in the case that $F RF^*$ is invertible reduces to
\begin{equation}
\label{eq:hatb2i}
\widehat{\beta}=
\frac{\tr((FRF^*)^{1/2})}{\tr(R)}.
\end{equation}  

The C-LSGUF vectors are then the columns of  $\widehat{\Phi}$ given by
\begin{equation}
\widehat{\Phi}=\hat{\beta}UV^* \Sigma \F^*,
\end{equation}  
where $\hat{\beta}$ is given by (\ref{eq:hatb2}). If $F RF^*$ is
invertible, then
\begin{equation}
\widehat{\Phi}=\hat{\beta}(FRF^*)^{-1/2}FR,
\end{equation}  
where $\hat{\beta}$ is given by (\ref{eq:hatb2i}).

\newpage
We summarize our results regarding constrained optimal GU frames in the
following theorem: 
\begin{theorem}[Constrained least-squares GU frames]
\label{thm:lsguf}
Let $\{ \varphi_i \}$ be a set of $n$ vectors in an $m$-dimensional
Hilbert space $\HH$ that span $\HH$, and let $F$ be the matrix
of columns $\varphi_i$.  Let $\{ \hat{\phi}_i \}$ denote the optimal
$n$ GU frame vectors that minimize the least-squares error defined by
(\ref{eq:serror})-(\ref{eq:error}), subject to the constraint
(\ref{eq:guconst}), and let $\widehat{\Phi}$ be the matrix of columns
$\phi_i$. 
Let $R$ be the matrix with $i$th  row equal to 
$P_i\{a_j,1 \leq j \leq n\}$, so that
$\beta^2R$ is the Gram matrix of inner products
$\{\inner{\phi_i}{\phi_j}\}$,
and let $\F$ be the Fourier transform matrix that
diagonalizes $R$.  Let $A$ be the
diagonal matrix with diagonal elements
$\{\alpha_j=n^{1/2}\hat{a}_j,1 \leq j \leq n\}$ where 
$\{\hat{a}_j,1 \leq j \leq n\}$
is the FT  of the
sequence $\{a_j,1 \leq j \leq n\}$, let $\Sigma$ be an $m \times n$
diagonal matrix with diagonal elements $\sqrt{\alpha_i}$ for values of
$i$ for which $\alpha_i \neq 0$, and let 
 $U$ and $V$ be the right-hand unitary
matrix and left-hand 
unitary matrix respectively in the SVD of $F\F \Sigma^*$.
Then,
\begin{enumerate}
\item if $\beta=\beta_0$ is given, then $\widehat{\Phi}=\beta_0 UV^*
\Sigma \F^*$. If in addition $FRF^*$ is invertible, then
$\widehat{\Phi}=\beta_0 (F\F A\F^*F^*)^{-1/2}F\F
A\F^*=\beta_0 (FRF^*)^{-1/2}F R$;
\item if $\beta>0$ is chosen to
minimize $E$, then 
$\widehat{\Phi}=\hat{\beta} UV^* \Sigma \F^*$, where $\hat{\beta}$ is given by
(\ref{eq:hatb2}). If in addition  $FRF^*$ is invertible, then
$\widehat{\Phi}=\hat{\beta} (F\F A\F^*F^*)^{-1/2}F\F
A\F^*=\hat{\beta}(F RF^*)^{-1/2}FR$, where $\hat{\beta}$ is given by
(\ref{eq:hatb2i}). 
\end{enumerate}
\end{theorem}  

\newpage
%%%%%%%%%%%%%%%%%%%%%%%%%%%%%%%%%%%%%%%%
\subsection{Least-squares GU frame}
\label{sec:gu}
%%%%%%%%%%%%%%%%%%%%%%%%%%%%%%%%%%%%%%

We now consider the least-squares problem  in which both the scaling
factor $\beta$ and
 the numbers
$\{a_j,1 \leq j \leq n\}$ in (\ref{eq:guconst}) are chosen to minimize $E$.
Thus, we seek a set of  vectors
$\{\phi_i\}$ that minimize the least-squares error $E$ of
(\ref{eq:serror})-(\ref{eq:error}) 
subject to 
\begin{equation}
\label{eq:guconst2}
\{\inner{\phi_i}{\phi_j},1 \leq j \leq n\}=P_i\{a_j,1 \leq j \leq n\},
\end{equation}
where $P_i\{a_j,1 \leq j \leq n\}$ is a known permutation of the
unknown numbers $\{a_j,1 \leq j \leq n\}$, chosen such that
the matrix $R$ whose $i$th row is equal to $P_i\{a_j,1 \leq j
\leq n\}$ is Hermitian, non-negative definite and diagonalized by a FT
matrix $\F$.

This problem has also been  considered in the context of
general least-squares inner product shaping \cite{E01}. It can be
shown  that the solution involves 
solving a problem of the form
\begin{equation}
\label{eq:guxy}
\max \,\sum_{i=1}^m |\inner{x_i}{y_i}|^2,
\end{equation} 
subject to
\begin{equation}
\label{eq:consty}
\inner{y_i}{y_j}=\delta_{ij},
\end{equation} 
where the vectors $\{x_i\}$ are known and are a function of the given
vectors $\{\varphi_i\}$.

As we now show, this  problem is equivalent to a quantum detection
problem, for which 
there is no known analytical solution in the general case.

%%%%%%%%%%%%%%%%%%%%%%%%%%%%%%%%%%
\subsubsection{Connection with quantum detection}
\label{sec:qd}
%%%%%%%%%%%%%%%%%%%%%%%%%%%%%%%%%

In a quantum detection problem,  a system is
prepared in one of $m$ known (pure) states that are described by
vectors $s_i$ in a  Hilbert space $\HH$, and the problem is
to  detect the state  prepared by performing a measurement on the
system. The measurement is described in terms  of a set of orthogonal
measurement 
vectors $q_i$.
Given a set of measurement vectors $q_i$, and assuming equal prior
probabilities on the different states, the 
probability of detection is given by \cite{H76},
\begin{equation}
\label{eq:pqd}
P_{qd}=\frac{1}{m}\sum_{i=1}^m|\inner{q_i}{s_i}|^2.
\end{equation}
Comparing (\ref{eq:pqd}) with (\ref{eq:guxy}) we see that finding a
set of orthogonal measurement vectors to maximize the probability of
detection is equivalent to the maximization problem of
(\ref{eq:guxy})-(\ref{eq:consty}).

Necessary and sufficient conditions for an optimum measurement
maximizing (\ref{eq:pqd}) have been derived \cite{H73,YKL75,H76}.
However, except in some particular cases \cite{H76,CBH89,BKMO97},
obtaining a closed-form analytical expression for the optimal
measurement directly from these conditions is a difficult and unsolved
problem.
Iterative algorithms
for maximizing $(\ref{eq:pqd})$ for arbitrary vector sets are given in
\cite{H82,E01}.  

We  conclude that in general there is no known analytical
expression for the GULSF. In practice, the GULSF may be obtained using 
the iterative algorithms of \cite{H82,E01}.  
A more detailed discussion on the GULSF can be found in \cite{E01}.

%%%%%%%%%%%%%%%%%%%%%%%%%%%%%%%%%%%%%%%%%%%%%%%%%%
\section{Distance Properties of GU Frames}
\label{sec:distance}
%%%%%%%%%%%%%%%%%%%%%%%%%%%%%%%%%%%%%%%%%%%%%%%%%

So far we have mainly been concerned with structural properties of GU
frames. In this section we study the Euclidean distance
properties of GU frames. 

Suppose we are given a GU frame $\{\phi_i=U_i \phi,1 \leq i \leq n\}$
generated by the group $\Q$
with $||\phi||^2=1$.
We would like to characterize the distance profile
$\alpha(i,j)=||\phi_i-\phi_j||^2$ 
for all $i,j$. 

Since the vectors $\phi_i$ are geometrically uniform, $\{\alpha(i,j),1 \leq
j \leq n\}$ is just a permutation of $\{d(i)=||\phi-\phi_i||^2,1
\leq i \leq n\}$. Furthermore,
\begin{equation}
d(i)=\inner{\phi}{\phi}+\inner{U_i\phi}{U_i\phi}-2\Re(\inner{U_i\phi}{\phi})=
2(1-\Re(a_i)),
\end{equation}
where $\{a_i,1 \leq i \leq n\}$ are the elements of the first row of the
Gram matrix corresponding to the frame $\{\phi_i\}$.
%These elements are determined
%by the FT matrix $\F$ and the eigenvalues $\lambda_i$ of $G$.

In  applications it may be desirable to construct a GU frame such
that $d(i)>0$ for $2 \leq i \leq n$. Since
$d(i)=||(I-U_i)\phi||^2$,  
a sufficient condition is that
\begin{equation}
\label{eq:fpf}
\mbox{det}(I-U_i) \neq 0, \quad 2 \leq i \leq n,
\end{equation}
which is satisfied if and only if none of  the  matrices
$U_i$ has an eigenvalue equal to $1$.  
Note that if (\ref{eq:fpf}) is satisfied, then $d(i)>0,2 \leq i \leq n$
regardless of the generating vector $\phi \,\neq\,0$.
Groups with unitary representations satisfying (\ref{eq:fpf}) are known
as {\it fixed-point free groups}, and have been
studied  extensively in the literature (see \eg \cite{JL93}). Thus if
$\Q$ is a representation of a fixed-point free group, then we have that
$d(i)>0,2 \leq i \leq n$.

Fixed-point free groups have recently been studied in the context 
of unitary space-time codes \cite{S00}. In particular, it was shown in \cite{S00}
that an abelian group of matrices $\{U_i\}$ satisfies (\ref{eq:fpf}) if
and only if it is 
cyclic, \ie $U_i=U^i$ with $U^n=I$, and where $U$ 
can be parameterized as
\begin{equation}
U=\mbox{diag}(e^{j2\pi u_{1}/n},\ldots,e^{j2\pi u_{n}/n}),
 \label{udef}
\end{equation}
where $u_{k}$ is relatively prime to $n$ for all $k$.  An
optimization over the $u_{k}$ can be performed to obtain distance
profiles with certain prescribed properties  (there are 
$\phi(n)$ positive integers less than $n$ that are relatively prime to
$n$, where $\phi(n)$ denotes the Euler totient function of $n$). 

\newpage

%%%%%%%%%%%%%%%%%%%%%%%%%%%%%%%%%%%%%%%%%%%%%%%%%%
\section{Conclusion}
%%%%%%%%%%%%%%%%%%%%%%%%%%%%%%%%%%%%%%%%%%%%%%%%%

In this paper we introduced the concept of GU and CGU frames and
discussed some of their key properties. These frames may be viewed as
generalizations of WH frames and filter bank frames to the group case.
A fundamental characteristic of these frames is that they
posses  strong symmetry properties that may be desirable in a variety of
applications. In particular, like WH frames and wavelet frames, GU
frames are generated by a single generating vector. Furthermore, the
canonical and dual frame vectors associated with a GU frame are themselves
GU and are therefore also generated by a single generating vector which
can be computed very efficiently using a Fourier transform matrix
defined over an appropriate group.

We also showed that GU frame vectors posses interesting symmetry
properties when one or more frame elements are removed. This property
of GU frames may be  of importance in multiple description
source-coding where it is often desirable  that the
quality of the reconstruction should not depend on the particular element
lost (removed).

Although in this paper we have focused on the
case in which the underlying group is a finite abelian group,
many of the results can be extended to the more general case of
infinite-dimensional  and non-abelian groups. 
An interesting direction for further research is to characterize these
more general cases of GU frames using 
possibly continuous-time Fourier transforms defined over non-abelian
groups (see 
\eg \cite{T99}).

%%%%%%%%%%%%%%%%%%%%%%%%%%%%%%%%
\section*{Acknowledgments}
%%%%%%%%%%%%%%%%%%%%%%%%%%%%%%%%

The authors wish to thank Prof.\ G.\ D.\ Forney for many fruitful
discussions.

\newpage
\begin{singlespace}
\bibliography{frames2}
\bibliographystyle{IEEEbib.bst}
\end{singlespace}

\end{document}